\documentclass[11pt]{article}
\usepackage[a4paper, margin=2cm]{geometry}

\usepackage[utf8]{inputenc}
\usepackage[T1]{fontenc}

\usepackage{bm}
\usepackage{physics}

\usepackage{tikz-cd}

\usepackage{upgreek}

\usepackage{amsmath,amsfonts,amssymb}

\usepackage{graphicx}
\usepackage{subfigure}

\usepackage{authblk}

\usepackage{hyperref}
\hypersetup{colorlinks,hidelinks}

\usepackage[backend=biber, style=numeric, giveninits=true, sorting=none]{biblatex}

\newcommand{\bigleft}{\bigl}
\newcommand{\bigright}{\bigr}

\newcommand{\leftf}{\mathopen{}\mathclose\bgroup\left}
\newcommand{\rightf}{\aftergroup\egroup\right}

\newcommand{\bigleftf}{\mathopen{}\mathclose\bgroup\bigl}
\newcommand{\bigrightf}{\aftergroup\egroup\bigr}

\newcommand{\ctilde}{\tilde}

\newcommand{\cbar}{\bar}

\newcommand{\ve}[1]{\mathbf{\bm{#1}}}
\newcommand{\mx}[1]{\mathbf{\bm{#1}}}

\newcommand{\yv}{\ve{y}}
\newcommand{\yvj}{\ve{y}^j}
\newcommand{\yvnull}{\ve{y_0}}
\newcommand{\yvjp}{\ve{y}^{j+1}}
\newcommand{\tj}{t^j}

\newcommand{\yvdot}{\dot{\ve{y}}}
\newcommand{\fv}{\ve{f}}
\newcommand{\RR}{\mathbb{R}}
\newcommand{\Dt}{\Delta t}

\newcommand{\h}{h}

\newcommand{\fvt}{\ve{\ctilde{f}}}
\newcommand{\yvt}{\ve{\ctilde{y}}}
\newcommand{\yvtdot}{\dot{\ve{\ctilde{y}}}}

\newcommand{\yvb}{\ve{\cbar{y}}}
\newcommand{\yvbdot}{\dot{\ve{\cbar{y}}}}
\newcommand{\fvb}{\ve{\cbar{f}}}
\newcommand{\yvbnull}{\ve{\cbar{y}_0}}

\newcommand{\Phihf}{\Phi_{\h,\fv}}
\newcommand{\Phihfb}{\Phi_{\h,\fvb}}

\newcommand{\PhihA}{\Phi_{\h,A}}
\newcommand{\PhihB}{\Phi_{\h,B}}
\newcommand{\PhihAinv}{\Phi_{\h,A}^{-1}}

\newcommand{\fvbt}{\ve{\ctilde{\cbar{f}}}}
\newcommand{\fbt}{{\ctilde{\cbar{f}}}}

\newcommand{\ft}{{\ctilde{f}}}
\newcommand{\yt}{{\ctilde{y}}}

\newcommand{\yb}{{\cbar{y}}}

\newcommand{\vtrf}{{\Psi}}
\newcommand{\vtrfinv}{{\Psi}^{-1}}

\newcommand{\D}{\mx{D}}

\newcommand\Ordo[1]{\mathcal{O}\leftf(#1\rightf)}

\newcommand{\qv}{\ve{q}}
\renewcommand{\pv}{\ve{p}}
\newcommand{\pvj}{\ve{p}^j}
\newcommand{\pvjp}{\ve{p}^{j+1}}
\newcommand{\pvjph}{\ve{p}^{j+1/2}}
\newcommand{\qvj}{\ve{q}^j}
\newcommand{\qvjp}{\ve{q}^{j+1}}

\newcommand{\zv}{\ve{z}}
\newcommand{\qvdot}{\dot{\ve{q}}}
\newcommand{\pvdot}{\dot{\ve{p}}}
\newcommand{\zvdot}{\dot{\ve{z}}}
\newcommand{\Hq}{H_q}
\newcommand{\Hp}{H_p}
\newcommand{\Hbar}{\cbar{H}}
\newcommand{\Htilde}{\ctilde{H}}
\newcommand{\Hbartilde}{\ctilde{\Hbar}}

\newcommand{\qvb}{\ve{\cbar{q}}}
\newcommand{\pvb}{\ve{\cbar{p}}}
\newcommand{\qvbj}{\ve{\cbar{q}}^{j}}
\newcommand{\pvbj}{\ve{\cbar{p}}^{j}}
\newcommand{\qvbjp}{\ve{\cbar{q}}^{j+1}}
\newcommand{\pvbjp}{\ve{\cbar{p}}^{j+1}}
\newcommand{\qvbjm}{\ve{\cbar{q}}^{j-1}}
\newcommand{\pvbjm}{\ve{\cbar{p}}^{j-1}}

\newcommand{\Qv}{\ve{Q}}
\newcommand{\Pv}{\ve{P}}

\newcommand{\Qvinv}{\ve{Q^{-1}}}

\newcommand{\Jcal}{\mathcal{J}}
\newcommand{\Jcalinv}{\mathcal{J}^{-1}}

\newcommand{\T}{\mathrm{T}}

\newcommand{\qb}{\cbar{q}}
\renewcommand{\pb}{\cbar{p}}

\newcommand{\Qinv}{Q^{-1}}

\newcommand{\al}{\alpha}
\newcommand{\be}{\beta}

\newcommand\Exp{\operatorname{Exp}}
\newcommand\Expof[1]{\Exp\leftf(#1\rightf)}
\newcommand\Id{\mathrm{Id}}

\newcommand{\Nm}{\mx{0}}
\newcommand{\Nv}{\ve{0}}
\newcommand{\Idm}{\mx{1}}

\newcommand{\Trp}{\mathrm{T}}

\newcommand{\Mminv}{\mx{M}^{-1}}

\newcommand{\Umod}{\hat{U}}

\newcommand{\Uqmhalf}{U_q^{-1/2}}
\newcommand{\Uqphalf}{U_q^{1/2}}

\newcommand{\epsxy}{\varepsilon_{(x,y)}}
\newcommand{\epsrphi}{\varepsilon_{(r,\varphi)}}

\newcommand{\Cm}{\mx{C}}

\newcommand{\Hcal}{H}
\newcommand{\Tcal}{T}
\newcommand{\Ucal}{U}
\newcommand{\Lcal}{L}
\newcommand{\pr}{p_{r}}
\newcommand{\pth}{p_{\theta}}

\newcommand{\px}{p_{x}}
\newcommand{\py}{p_{y}}

\newcommand\atantwo{\operatorname{atan2}}

\newcommand\deltaQ[2]{\Xi_{#1, #2}}
\newcommand\deltaHpHq{\deltaQ{H_p H_q}{\Qv}}
\newcommand\deltabarHpHq{\deltaQ{\Hbar_{\pb} \Hbar_{\qb}}{\Qvinv}}

\begin{document}

\title{On the coordinate system-dependence of the accuracy of symplectic numerical methods}

\author[1,2]{Don\'{a}t M. Tak\'{a}cs\thanks{takacs@energia.bme.hu}}
\author[1,2]{Tam\'{a}s F\"{u}l\"{o}p}
\affil[1]{Department of Energy Engineering, Faculty of Mechanical Engineering, Budapest University of Technology and Economics, M\H{u}egyetem rkp.~3, Budapest, H-1111, Hungary}
\affil[2]{Montavid Thermodynamic Research Group, Society for the Unity of Science and Technology, Lovas \'{u}t~18., Budapest, H-1012, Hungary}
 

\maketitle

\begin{abstract}
Symplectic numerical methods have become a widely-used choice for the accurate simulation of Hamiltonian systems in various fields, including celestial mechanics, molecular dynamics and robotics. Even though their characteristics are well-understood mathematically, relatively little attention has been paid in general to the practical aspect of how the choice of coordinates affects the accuracy of the numerical results, even though the consequences can be computationally significant.
The present article aims to fill this gap by giving a systematic overview of how coordinate transformations can influence the results of simulations performed using symplectic methods. We give a derivation for the non-invariance of the modified Hamiltonian of symplectic methods under coordinate transformations, as well as a sufficient condition for the non-preservation of a first integral corresponding to a cyclic coordinate for the symplectic Euler method. We also consider the possibility of finding order-compensating coordinate transformations that improve the order of accuracy of a numerical method. Various numerical examples are presented throughout.
\end{abstract}



\section{Introduction}  
Symplectic numerical methods have become a dependable approach for the simulation of Hamiltonian systems in cases where accurate and qualitatively correct numerical solutions are needed, including very long integration times. Notable applications of symplectic methods include
precise calculation of the orbits of celestial objects and spacecraft \cite{duncan1995dynamical,chyba2008role,laskar2011la,zeebe2023orbitn},
solving optimal control problems in robotics \cite{pekarek2007discrete,oberblobaum2011discrete,shareef2016simultaneous,wang2022iterative},
determining geological timescales for climate models \cite{zeebe2019solar,zeebe2022geologically}, simulation of biological systems \cite{beck2003symplectic,farago2018operator} and
calculating the dynamic aperture of particle accelerators \cite{forest2006geometric,ruth1983canonical,forest1989canonical,erdelyi2001optimal},
among others. More recently, symplectic methods and their extensions have also been applied successfully to stochastic systems \cite{chen2021asymptotically,dambrosio2024strong} and for the time integration of partial differential equations with various spatial discretizations \cite{morrison2017structure,uzunca2023global,takacs2024thermodynamically,jones2024discrete}. Behind the success of symplectic methods lies the mathematical fact that, unlike most other numerical methods, they have a corresponding Hamiltonian system which they integrate exactly \cite{benettin1994hamiltonian,hairer2006geometric}. This so-called modified or distorted\footnote{Following the reasoning given in \cite{takacs2024improving}, in this article we use the term \emph{distorted} instead of the more commonly used \emph{modified}, when referring to expressions obtained by backward error analysis of a numerical method.} Hamiltonian is close to the Hamiltonian of the original system to the order of $\Ordo{h^r}$, where $r$ is the order of convergence of the method. In the past decades, significant mathematical effort has been put into the deep understanding of the underlying distorted Hamiltonians of symplectic methods and their subsequent extensions \cite{forest2006geometric,leimkuhler2005simulating,hairer2006geometric,feng2010symplectic,feng1991note,kane2000variational,moan2006modified,morrison2017structure}.

At the same time, processing for symplectic methods has been developed, which uses canonical transformations close to the identity to perform the numerical time stepping in phase space coordinates that enables a higher order of accuracy \cite{takahashi1984monte,rowlands1991numerical,lopezmarcos1997explicit,mclahlan2002splitting,casas2010processed}. Despite the success of these methods, to the best of our knowledge, there has been sparse work investigating how the choice of coordinates affects the performance of symplectic methods.

The appropriate choice of coordinates in the formulation of Hamiltonian systems has a long history, since a proper choice can enable the solution of problems not solvable in the original formulation, as first shown by Jacobi \cite{jacobi1836uber,nakane2002early,arnold1989mathematical,goldstein2001classical}. Due to this, the question of what coordinates to use in the formulation of Hamiltonian equations is primarily viewed as one of convenience: after all, the Hamiltonian describing the system is invariant under canonical transformations by definition. This is especially true for the choice of generalised coordinates%
\footnote{
    In the remaining part of this article, for the sake of conciseness and unambiguity, we will refer to canonical transformations induced by
a point transformation of the generalised coordinates
as \emph{coordinate transformations}, and refer to the general transformations of phase space coordinates as \emph{canonical transformations}. The term \emph{variable transformations} is reserved for ODEs that are not necessarily Hamiltonian.
}
(which induce a point transformation of the phase space), which can be an ad hoc or intuitive choice made early during the formulation of the problem. However, the improved performance of processing methods shows that symplectic numerical methods are not indifferent to the coordinates used in the formulation of the original problem.

The main goal of this paper is to give a systematic treatment of how the choice of (generalised) coordinates and the induced coordinate transformations affect the accuracy of symplectic methods, and introduce novel results that help in distinguishing between coordinate transformations that have an effect on the accuracy of results and those that do not. In Section~\ref{sec:vtrf}, we introduce the necessary formalism and show how variable transformations in general can affect the numerical solutions of ODEs. In Section~\ref{sec:processingsym}, we give a brief review of processing methods for symplectic methods. Then in the remaining part of the article, we investigate two main aspects of symplectic methods under coordinate transformations: the non-invariance of the distorted Hamiltonian and the preservation of first integrals. To the best of our knowledge, the present treatment of these questions is a novel approach not encountered in the relevant literature. For the distorted Hamiltonians of symplectic methods, we show in Section~\ref{sec:distortedhamiltonian} that they are not invariant to coordinate transformations in general, and demonstrate this result numerically, as well as pointing out possible approaches for exploiting this fact to achieve better accuracy. Regarding first integrals, in Section~\ref{sec:firstintegrals}, we derive and demonstrate a sufficient condition for the non-preservation of a first integral corresponding to a cyclic coordinate for the symplectic Euler method.



\section{Variable transformations in numerical methods for ODEs}\label{sec:vtrf}
\subsection{On variable transformations in general}

Let us consider an $n$-dimen\-sional, autonomous ordinary differential equation with an initial value, written as
\begin{gather}
    \yvdot = \fv \leftf( \yv \rightf);\quad \yv\leftf(0\rightf) = \yvnull, \label{eq:ode}
\end{gather}
having $ \yv : \RR \rightarrow \RR^n $ as the solution with initial value $\yvnull$, and a generating vector field $\fv : \RR^n \rightarrow \RR^n $, where the overdot denotes differentiation with respect to time. Consider a (sufficiently smooth, one-to-one) variable transformation $\vtrf : \RR^n \rightarrow \RR^n$, which transforms $\yv$ to
\begin{gather}
    \yvb = \vtrf\leftf( \yv \rightf).\label{eq:vartrf}
\end{gather}

Inserting \eqref{eq:vartrf} into \eqref{eq:ode} and rearranging yields the transformed problem
\begin{gather}
    \yvbdot = \left( \D \vtrfinv \leftf( \yvb \rightf) \right)^{-1} \fv \circ \vtrfinv \leftf( \yvb \rightf) =: \fvb(\yvb); \quad \yvb(0) = \vtrf(\yvnull) =: \yvbnull, \label{eq:odetrf}
\end{gather}
where $\D$ denotes the derivative.

A one-step numerical method $\Phihf: \RR^n \rightarrow \RR^n$ with time step $\h$ generates the approximate solution to \eqref{eq:ode} as
\begin{gather}
    \yvjp = \Phihf \leftf( \yvj \rightf)
\end{gather}
where the approximate solution is denoted as $\yvj$ at time instant $\tj := j \h $. Accordingly, a numerical method is invariant with respect to the variable transformation if and only if
\begin{gather}
    \Phihfb \equiv \vtrf \circ \Phihf \circ \vtrfinv \label{eq:vtrfinvariant}
\end{gather}
is satisfied. 

Usually, \eqref{eq:vtrfinvariant} is known to hold for certain classes of $\vtrf$ for most methods, such as for linear transformations in the case of Runge--Kutta methods \cite{hairer2002solving2}, and is sometimes exploited during the analysis of numerical methods \cite{lubich1990convergence}.

However, in what follows, we are considering transformations for which \eqref{eq:vtrfinvariant} does \emph{not} hold, in such a way that it changes the accuracy of $\Phihfb$ in a certain sense. One of the main ways to achieve this is to find a $\vtrf$ such that it increases the order of convergence of the numerical integration: this is referred to here in general as \emph{order compensation} \cite{takacs2024improving}. Transformations that achieve order compensation can be classified according to whether they depend on $\h$ or not: in the following subsection, we consider the latter, and will show examples of the former in Section~\ref{sec:processingsym}.

\subsection{Order compensation using variable transformation}

Using backward error analysis \cite{griffiths1986scope,lopezmarcos1996cheap,reich1999backward,hairer2000asymptotic,moan2006modified}, it can be shown that for a discrete-time one-step method $\Phihf$, a corresponding continuous-time system can be constructed for which $\Phihf$ is an exact integrator. This system is described by the so-called modified or distorted equation
\begin{gather}
    \yvtdot = \fvt\leftf( \yvt \rightf); \quad \yvt(0) = \yvnull
\end{gather}
where, by definition, the distorted solution fulfils the condition
\begin{gather}
    \yvt\leftf(\tj\rightf) = \yvj,
\end{gather}
and the distorted vector field (DVF) $\fvt$ has the form of an (actually asymptotic) power series
\begin{gather}
    \fvt = \fv + \h \;\! \fvt_1 + \h^2 \;\! \fvt_2 + \ldots
\end{gather}
for a consistent method $\Phihf$ giving the exact solution in the $h \rightarrow 0$ limit. The terms $\fvt_i$ depend both on the original vector field and on the numerical method used.

Consequently, if we use a suitable variable transformation $\Psi$ (with regards to which the numerical solution is not invariant), one or more terms in the DVF $\fvbt$ of the transformed system $\fvb$ could be cancelled. Using such an approach, the order of accuracy of the original method can be improved. For example, a compensation of second order would require $\fvbt_1 = \Nm$.

Furthermore, since the transformation only has to be used for the initial condition and (in the inverse direction) at time steps where the actual result is needed, the numerical calculation as a whole incurs negligible additional cost. However, determining the suitable transformation itself is often not trivial, as we will show in the following example.

\subsection{Two examples: explicit Euler method, first-order system}
Let us use the explicit Euler method
\begin{gather}
    \yvjp = \underbrace{\yvj + \h \;\! \fv \leftf( \yvj \rightf)}_{\Phihf(\yvj)},
\end{gather}
for solving \eqref{eq:ode} for $n=1$.

Then, according to backward error analysis~\cite{hairer2000asymptotic,takacs2024improving}, we have the DVF
\begin{gather}
    \ft(\yt) = f(\yt) - \frac{1}{2} \h \;\! f'(\yt) f(\yt) + \Ordo{\h^2}.
\end{gather}
Similarly, calculating the DVF of the explicit Euler method applied to the transformed system \eqref{eq:odetrf}, and transforming back to the original variables yields
\begin{gather}
    \bigleft(\vtrf'(\yt) \bigright)^{-1} \fbt \circ \vtrf(\yt) =
    f(\yt) - \frac{1}{2} h \left( f'(\yt) f(\yt) + \frac{\vtrf''(\yt)}{\vtrf'(\yt)} f^2(\yt) \right) + \Ordo{\h^2}. \label{eq:DVFbarunbar}
\end{gather}
Immediately, we can see that the distorted equations corresponding to the transformed and the original systems are equivalent if $\vtrf''\equiv0$, i.e., if $\vtrf$ is linear. (Here, this is only apparent for the first order of $\h$, but can be generalised to higher-order terms as well.) This means that in this case, the
necessary condition for the invariance \eqref{eq:vtrfinvariant} is the linearity of $\vtrf$. Otherwise, the discrepancy between the two DVFs opens up a possibility for order compensation. According to \eqref{eq:DVFbarunbar}, the condition for second-order compensation takes the form of an ordinary differential equation for $\vtrf$:
\begin{gather}
    f'(y) f(y) + \frac{\vtrf''(y)}{\vtrf'(y)} f^2(y) = 0, \label{eq:1dcomp}
\end{gather}
which has the solution
\begin{gather}
    \vtrf(y) = \int_{y_0}^{y} \frac{C_1}{f(u)} \, \mathrm{d}u + C_2 \label{eq:vtrfexpr}
\end{gather}
with integration constants $C_1$ and $C_2$, and where the initial value $y_0$ has also been used.

\subsubsection{Example 1: explicit Euler method, first-order system}
\begin{figure}[htb]
    \centering
    \includegraphics[width=0.8\textwidth]{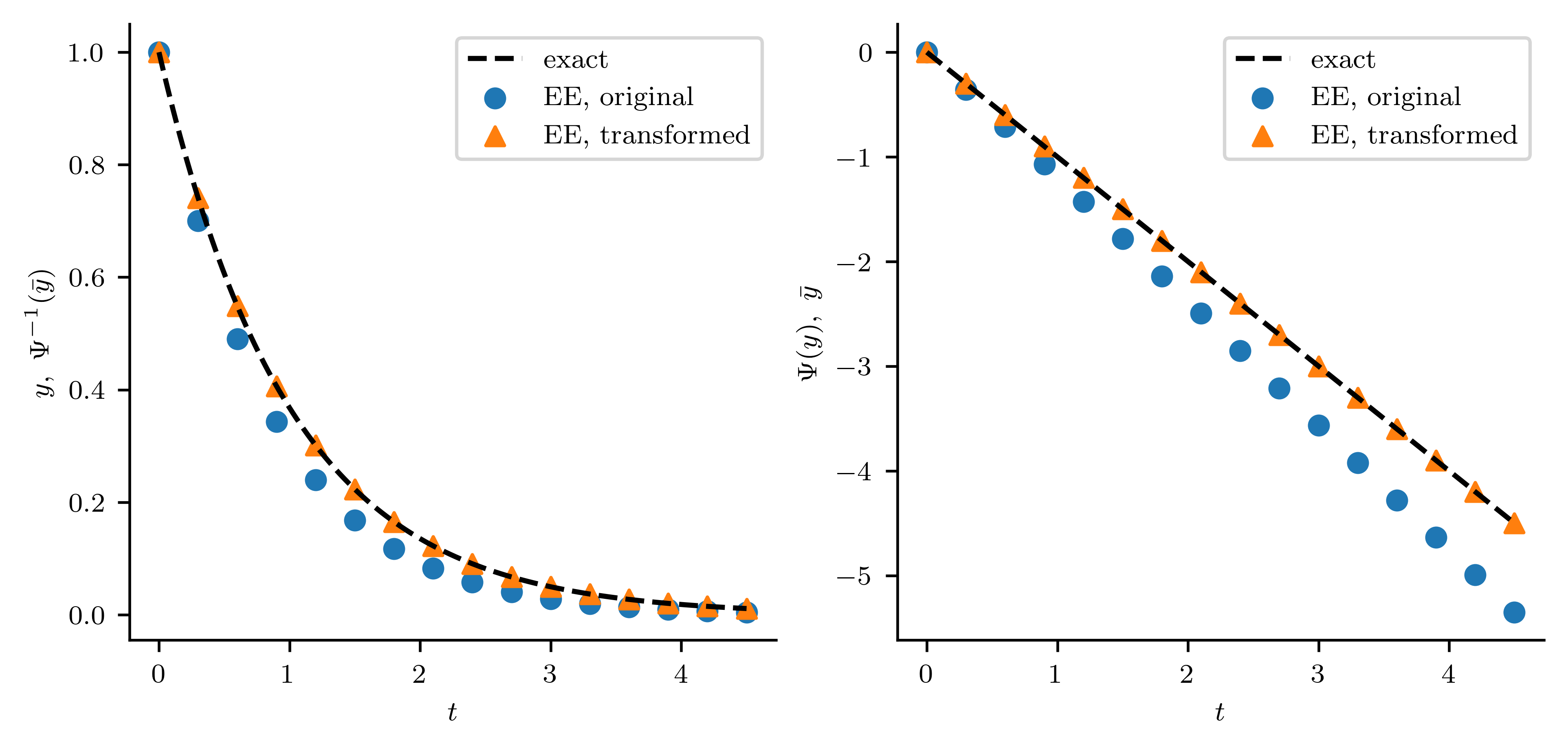}
    \vspace*{-4ex}
    \caption{Accuracy of the explicit Euler method for solving \eqref{eq:cooling} in the original and in the optimal coordinates, with $y_0=1$, $\alpha=1$ and $h=0.3$. Both numerical results and the exact result are shown in both coordinates.}
    \label{fig:exactcompEE}
\end{figure}

For certain vector fields, such as the ODE
\begin{gather}
    \dot{y} = - \alpha y \label{eq:cooling}
\end{gather}
with a positive constant $\alpha$ and initial condition $y_0 > 0$, expression \eqref{eq:vtrfexpr} can be evaluated to get a closed-form expression. In this case we obtain
\begin{gather}
    \vtrf(y) = \frac{1}{\alpha} \ln{y}
\end{gather}
as a suitable compensating variable transformation (with a proper selection of $C_1$ and $C_2$). Applying this variable transformation to \eqref{eq:cooling} gives the transformed ODE
\begin{gather}
    \dot{\yb} = - \alpha \label{eq:coolingtrf}
\end{gather}
which, integrated numerically by the explicit Euler method, actually gives not only a first-order accurate, but an exact solution to the original ODE. Fig.~\ref{fig:exactcompEE} compares the performance of the EE method solving \eqref{eq:cooling} in the original and the optimal coordinate system.

\subsubsection{Example 2: explicit Euler method, nonlinear system}
\begin{figure}[htb]
    \centering
    \includegraphics[width=0.8\textwidth]{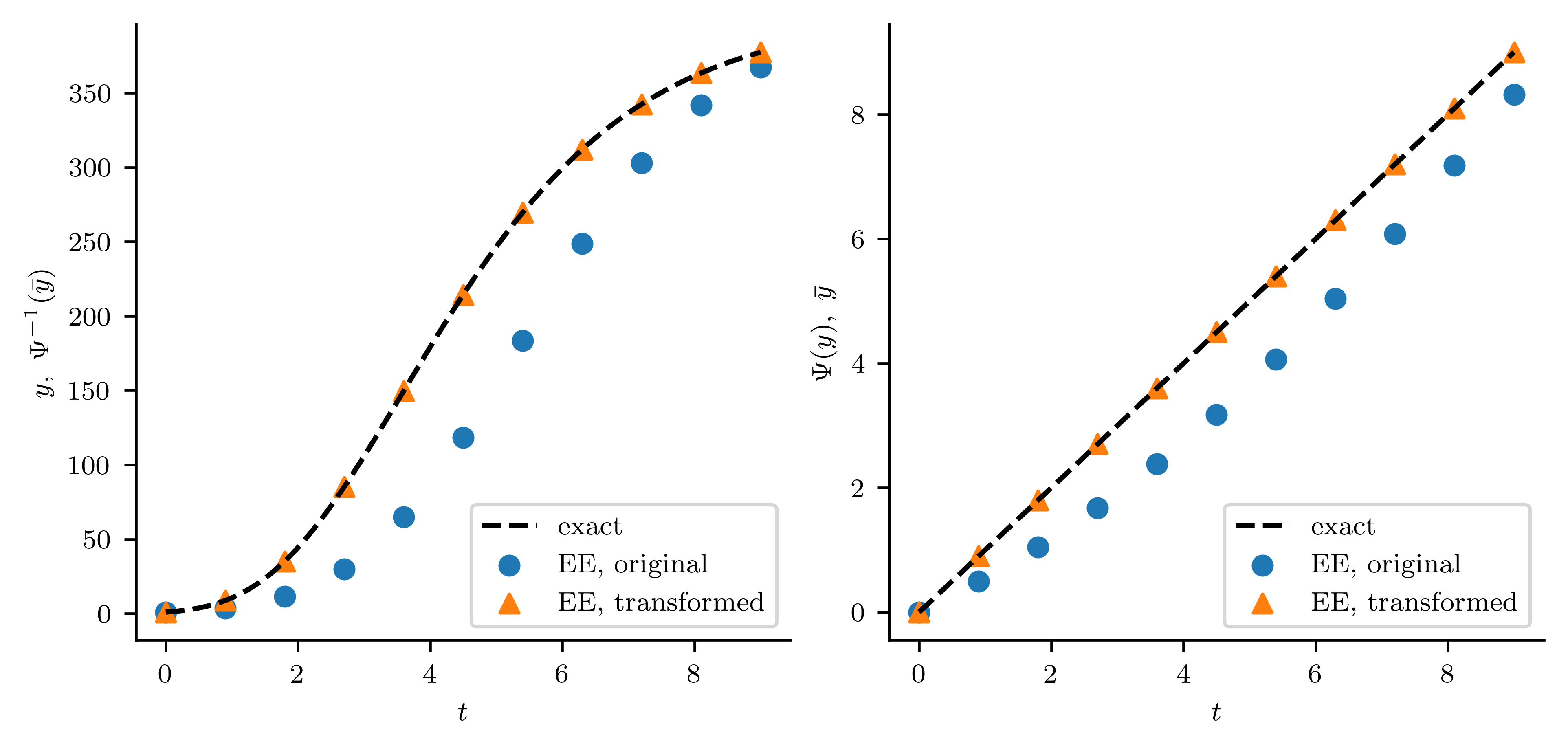}
    \vspace*{-4ex}
    \caption{Accuracy of the explicit Euler method for solving \eqref{eq:gompertz} in the original and in the optimal coordinates, with $y_0=3$, $a=2$, $b=0.5$, and $h=0.9$. Both numerical results and the exact result are shown in both coordinates.}
    \label{fig:exactcompEEgompertz}
\end{figure}

As a more complex example, consider the nonlinear Gompertz equation used in biology for describing the growth of plants, animals, bacteria and cancer cells \cite{gompertz1825nature,winsor1932gompertz,tjorve2017use,larsen2017models}, expressed as
\begin{gather}
    \dot{y} = y ( a - b \ln(y)), \label{eq:gompertz}
\end{gather}
with non-zero growth parameters $a$ and $b$ and initial condition $y_0 > 0$. Using again the explicit Euler method to solve it, evaluating \eqref{eq:vtrfexpr} yields
\begin{gather}
    \vtrf(y) = \frac{1}{b} \left( \ln\leftf( 1 - \frac{b}{a} \ln(y_0) \rightf) - \ln\leftf( 1 - \frac{b}{a} \ln(y) \rightf) \right)
\end{gather}
as an optimal transformation (with $C_1=1$ and $C_2=0$),
which can be shown to be invertible for the range of the solution. Consequently, the transformed ODE becomes the trivial $\dot{\yb} = 1$,
which, again, can be integrated exactly using the explicit Euler method, giving an exact numerical solution through the use of the optimal coordinate system. Fig.~\ref{fig:exactcompEEgompertz} shows a comparison with a solution obtained in the original coordinate system that exhibits significant errors at the same timestep.

At this point, it must be noted that analytically solving the original differential equation \eqref{eq:cooling} entails the calculation of a virtually identical integral to that of \eqref{eq:vtrfexpr}, thus, in a practical sense, we have not gained a lot from the transformation in this general form. However, this does not mean that the general approach can not be viable for specific cases.

\subsection{Viability of order compensation}

So far, we have de\-mon\-stra\-ted that there is potentially a lot to be gained from a well-chosen variable transformation aimed at order compensation, and that such a transformation can be systematically identified using backward error analysis. However, we have also shown that determining this transformation even in a one-dimensional case has a difficulty comparable to the original problem itself.

For a higher-dimensional case of \eqref{eq:ode}, it can be shown that the condition \eqref{eq:1dcomp} becomes
\begin{gather}
    \D \fv \:\! \fv + \left( \D \vtrfinv \right) \D^2 \vtrf \leftf( \fv, \fv \rightf) \equiv \Nv,
    \label{eq:ndcomp}
\end{gather}
which is essentially a partial differential equation for $\vtrf$. This makes the task of finding a coordinate-system even more complicated, a task that might have, in general, a difficulty higher than that of the original problem.

Nevertheless, this difficulty can be relaxed by using a combination of two approaches: we consider only a subset of possible problems (such as Hamiltonian systems and symplectic methods), and we use $\h$-dependent transformations close to the identity that only cancel terms up to a certain order. The so-called \emph{processing} method for symplectic methods is one such concept.

\section{Processing for symplectic methods}\label{sec:processingsym}
\subsection{Short overview of symplectic methods}
%
Symplectic numerical methods \cite{devogelaere1956methods,forest2006geometric,leimkuhler2005simulating,hairer2006geometric,feng2010symplectic,nasab2021partitioned} are a special class of numerical methods suitable for integrating Hamiltonian systems, where the preservation of the symplectic structure characteristic of the original system is important in most applications. A numerical method that is symplectic preserves various key properties of Hamiltonian systems: the geometry of the phase space, closed orbits, while guaranteeing a near-preservation of the total energy of the system. On the other hand, non-symplectic methods usually exhibit numerical dissipation (e.g.~the fourth-order Runge--Kutta method) or anti-dissipation (e.g.~the explicit Euler method) \cite{hairer2006geometric,takacs2024improving}, which makes them unsuitable for long-term calculations.

Equations describing autonomous Hamiltonian systems are of the form
\begin{gather}
    \underbrace{
    \begin{pmatrix}
        \qvdot \\ \pvdot
    \end{pmatrix}
    }_{\zvdot}
    =
    \underbrace{
    \begin{pmatrix}
        \Nm & \Idm \\
        -\Idm & \Nm \\
    \end{pmatrix}
    }_{\Jcal}
    \underbrace{
    \begin{pmatrix}
        \Hq(\qv, \pv) \\
        \Hp(\qv, \pv) \\
    \end{pmatrix}
    }_{\D H(\zv)}
\end{gather}
where $H(\qv, \pv): \RR^n \cross \RR^n \rightarrow \RR$ is the Hamiltonian as a function of the generalised coordinates and momenta $\qv \in \RR^n$ and $\pv \in \RR^n$, with $H_q$ and $H_p$ being shorthands for partial derivatives of $H$. In the following, we only consider Hamiltonians independent of $t$, corresponding to autonomous systems.

A numerical method $\Phihf$ is symplectic if and only if \cite{hairer2006geometric}
\begin{gather}
    \left(
        \pdv{\Phihf}{\zv} 
    \right)^{\T}
    \Jcalinv
    \left(
        \pdv{\Phihf}{\zv} 
    \right)
    =
    \Jcalinv
\end{gather}
holds.

As symplectic methods usually require a fixed time step size, they yield themselves readily to compensations or corrections that are otherwise not usable for variable-step methods. \emph{Processing} methods have this requirement of fixed step size: correspondingly, even though the original idea dates back to Butcher \cite{butcher1969effective}, they gained renewed popularity as geometric (and as a special case, symplectic) integrators became more widely used.

\subsection{Processing methods}\label{sec:processing}
Processing methods are a special case of splitting methods \cite{mclahlan2002splitting}, being of the form
\begin{gather}
    \Phihf = \PhihAinv \circ \PhihB \circ \PhihA, \label{eq:processing}
\end{gather}
where the three stages correspond to the \emph{preprocessor}, \emph{kernel} and \emph{processor}, respectively \cite{lopezmarcos1996cheap}. (The pre- and postprocessors are sometimes also called \emph{correctors}, making the methods $\Phihf$ and $\PhihB$ conjugate to each other \cite{mclahlan2002splitting}.)

Naturally, such methods have the property that
\begin{gather}
    \Phihf^n \equiv \PhihAinv \circ \PhihB^n \circ \PhihA \quad (n = 1, 2, \ldots),
\end{gather}
meaning that only the kernel needs to be evaluated repeatedly. Meanwhile, the preprocessor only needs to be evaluated once for the initial condition, and the postprocessor only at those instants where the actual output is needed. What's more, there are certain situations where even the repeated application of the kernel is sufficient \cite{takahashi1984monte,casas2010processed} without the pre- or postprocessors, as the calculated quantity of interest is invariant with respect to processing. Due to the above, it is customary to assume that the costs of evaluating the pre- and postprocessor stages can be neglected compared to the cost of the repeated application of the kernel \cite{lopezmarcos1996cheap,mclahlan2002splitting}.

It is also useful to introduce the relevant concept of \emph{effective order} \cite{butcher1969effective} here: a method $\Phi_h$ is said to be of effective order $k$, if there exists a corrector $\chi$ such that $\chi^{-1} \circ \Phi \circ \chi$ is of order $k$, while the method $\Phi_h$ is originally of order $m < k$. In exponential form, this definition can be expressed as
\begin{gather}
    \Expof{h \:\! \fv + \Ordo{h^{p+1}}} = \Expof{- h A} \Expof{h \;\! \fv + \Ordo{h^{m+1}}} \Expof{h A},
\end{gather}
where $\Exp$ is the exponential map of vector fields.

A fairly known example of an integrator having higher effective order is the Takahashi--Imada integrator \cite{takahashi1984monte} used in molecular dynamics simulations. Additionally, some splitting methods also have an effective order higher than their original order, explaining some of their serendipitous properties \cite{mclahlan2002splitting}. As a closely related fact, it should be noted here that a processing method \eqref{eq:processing} is symplectic if both $A$ and $B$ are symplectic, since their composition will also be symplectic.

Based on the above definitions, an order-compensation method that uses processing is a method of the form \eqref{eq:processing} that yields a method of higher effective order than its original order, while the processing stages are also known and are feasible to compute. There exist such approaches for symplectic methods \cite{rowlands1991numerical,lopezmarcos1996cheap,wisdom1996symplectic,lopezmarcos1997explicit}, and an important characteristic of all their approaches is that the corrector is of the form $\chi = \Id + \Ordo{h}$. This is significantly restricted compared to a general variable transformation $\vtrf$ (as outlined in Section~\ref{sec:vtrf}), yet, useful results can be achieved.

\subsection{Example: Rowlands time stepping}\label{sec:example_rowlands}
As an illustration of the previously explored ideas, in this section we show the construction and performance of the Rowlands time stepping (first discovered by \cite{rowlands1991numerical}, then rediscovered by \cite{wisdom1996symplectic}, and subsequently improved by \cite{lopezmarcos1997explicit}), compared to the St\"{o}rmer--Verlet method. 

The Rowlands method is applicable to separable Hamiltonians of the form
\begin{gather}
    H(\qv, \pv) = \underbrace{\frac{1}{2} \pv^\Trp \Mminv \pv}_{T(\pv)} + U(\qv),\label{eq:rowlandscondition}
\end{gather}
for which the (second-order) St\"{o}rmer--Verlet method can be written as
\begin{align}
    \pvjph &= \pvj - \frac{h}{2} U_q\leftf( \qvj \rightf) \label{eq:stormerverlet1}\\
    \qvjp  &= \qvj + h \Mminv \pvjph \label{eq:stormerverlet2}\\
    \pvjp &= \pvjph - \frac{h}{2} U_q\leftf( \qvjp \rightf),\label{eq:stormerverlet3}
\end{align}
where $U_q$ is shorthand for partial derivatives of $U$. Or, in exponential form, it can also be written as
\begin{gather}
    \Expof{\frac{h}{2} U} \Expof{h T} \Expof{\frac{h}{2} U}.\label{eq:stormerverletexp}
\end{gather}

The Rowlands method introduces an effective potential $\Umod$, defined as
\begin{gather}
    \Umod =  U - \frac{h^2}{24} U_q \Mminv U_q,
\end{gather}
using which the Rowlands time-stepping is defined (analogously to \eqref{eq:stormerverletexp}) as
\begin{gather}
    \Expof{\frac{h}{2} \Umod} \Expof{h T} \Expof{\frac{h}{2} \Umod}.\label{eq:rowlandsexp}
\end{gather}
The main advantage of this approach is that it increases the effective order of the Rowlands method to four, from the second-order accuracy of the St\"{o}rmer--Verlet method. In other words, \eqref{eq:rowlandsexp} can be viewed as a kernel of a fourth-order processing method. The pre- and postprocessing transformations were originally given by Rowlands as
\begin{align}
    \qvb &= \qv + \frac{h^2}{24} \Mminv U_q\leftf( \qv \rightf), \\
    \pvb &= \pv - \frac{h^2}{24} U_{qq}\leftf( \qv \rightf) \Mminv \pv,
\end{align}
but later it has been shown \cite{lopezmarcos1996cheap,lopezmarcos1997explicit} that these transformations can be approximated by the following expressions that are cheaper to evaluate. For the preprocessing step, this is
\begin{align}
    \qvbj &= \qvj + \frac{h^2}{24} \Mminv \left( \Uqphalf + \Uqmhalf \right), \\
    \pvbj &= \pvj - \frac{h}{12}\left(  \Uqphalf - \Uqmhalf \right),
\end{align}
with $U_q^{\pm 1/2} = U_q\left( \qvj \pm \frac{h}{2} \Mminv \pvj \right)$, and for the postprocessing
\begin{align}
    \qvj &= \qvbj + \frac{1}{12}\left( \qvbjp - 2 \qvbj + \qvbjm \right), \\
    \pvj &= \pvbj - \frac{1}{12}\left( \pvbjp - 2 \pvbj + \pvbjm \right).
\end{align}

As an example, we demonstrate the convergence of the St\"{o}rmer--Verlet and Rowlands methods with the dynamics of the planar elastic pendulum, a nonlinear system with two degrees of freedom also used in modelling various engineering problems (e.g. in mitigating the vibrations of power lines \cite{farzaneh2008atmospheric,kollar2022ice} or in modelling the motion of ships \cite{lee1992domains,lee1997chaotic}). The details of the mechanical model
is given in \ref{app:elasticpendulum}.

\begin{figure}[htbp]
    \centering
    \includegraphics[width=\textwidth]{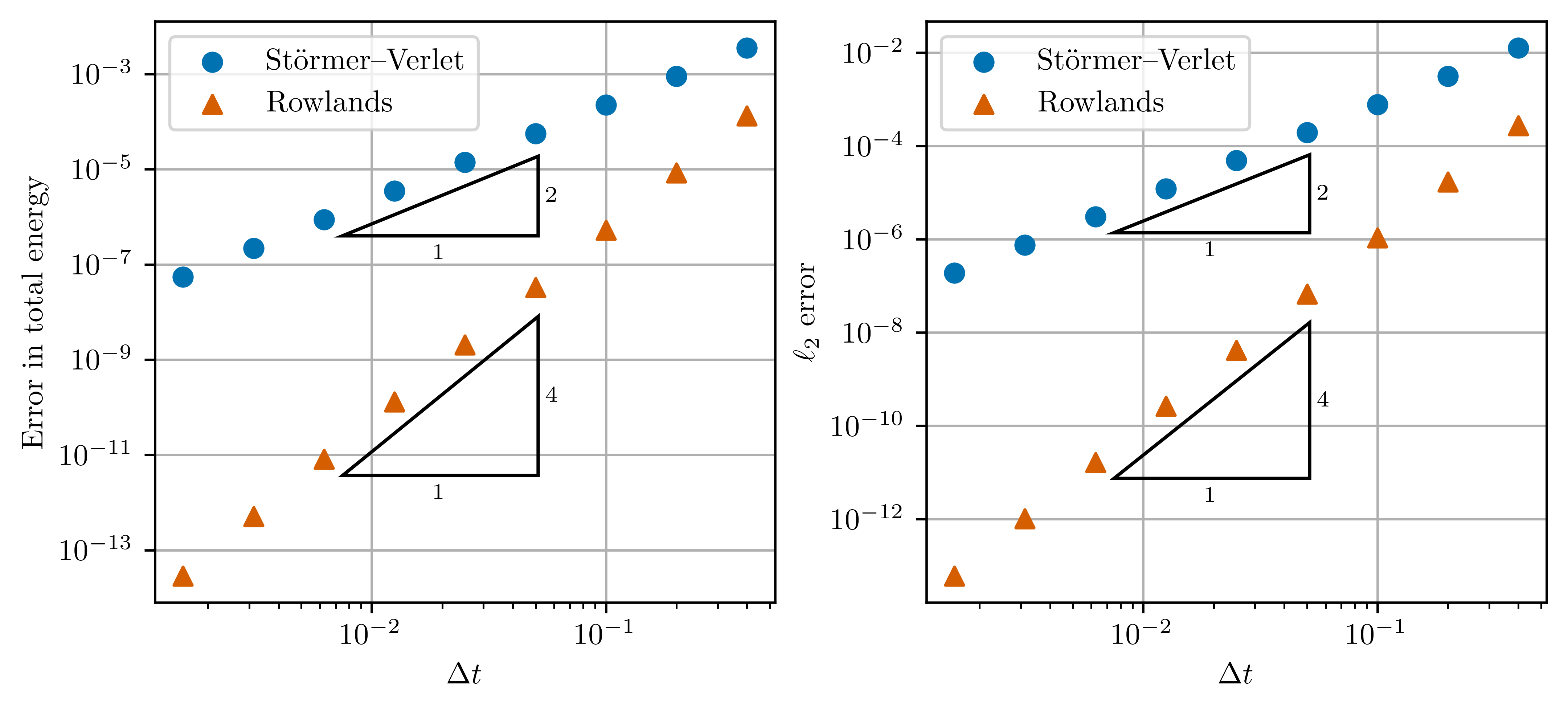}
    \vspace*{-4ex}
    \caption{Demonstration of the fourth-order convergence of the Rowlands method versus the second-order convergence of the St\"{o}rmer--Verlet method through the example of the planar elastic pendulum in Cartesian coordinates. (Parameters: $l{=}1,\, m{=}1,\, k{=}1,\, g{=}0.2,\, t_\mathrm{max}{=}4$.)}
    \label{fig:rowlandsconvergence}
\end{figure}

Fig.~\ref{fig:rowlandsconvergence} shows that the Rowlands method clearly outperforms the St\"{o}rmer--Verlet method both in terms of accuracy in the total energy and the $\ell_2$ error of the phase space coordinates, as it achieves a fourth-order convergence, compared to the latter method which is second-order. At the same time, the Rowlands method fundamentally uses the same time-stepping as the St\"{o}rmer--Verlet method, its overhead only being the one-time calculation of the effective potential $\Umod$, the one-time calculation of the preprocessing step, and the postprocessing at the where the results are to be evaluated.

\section{Effect of change of generalised coordinates on the distorted Hamiltonian}\label{sec:distortedhamiltonian}
In the previous section, we gave an overview of existing processing approaches for symplectic methods. These used processors that were $h$-flows of a Hamiltonian system close to the identity: in the following, we explore the more general class of transformations that are not close to the identity, and thus are not necessarily $h$-flows.

One of the key properties of Hamiltonian systems is that the Hamiltonian $H$ is invariant to canonical transformations \cite{goldstein2001classical}. For a time-independent canonical transformation
\begin{gather}
    \begin{pmatrix}
        \qvb \\ \pvb
    \end{pmatrix}
    =
    \begin{pmatrix}
        \Qv\leftf( \qv, \pv \rightf) \\
        \Pv\leftf( \qv, \pv \rightf)
    \end{pmatrix},
\end{gather}
the Hamiltonian in the transformed coordinates is
\begin{gather}
    \Hbar \bigleftf( \qvb, \pvb \bigrightf) \equiv H \bigleftf( 
        \Qv\leftf( \qv, \pv \rightf),
        \Pv\leftf( \qv, \pv \rightf)
    \bigrightf).
\end{gather}

At the same time, for symplectic methods, an important result from backward error analysis which enables the understanding of their properties is that symplectic methods are exact integrators of an underlying distorted system, which, additionally, is also described by a distorted Hamiltonian \cite{benettin1994hamiltonian}, denoted in the following as $\Htilde$. However, a crucial question not raised in the literature is whether a symplectic numerical method applied to the same system in two different coordinates yields equivalent distorted Hamiltonians. In other words: is the distorted Hamiltonian, resulting from discretization, invariant with respect to canonical transformations of the simulated system, i.e., is $\Htilde$ the same as $\Hbartilde$? As we will show in the following, the answer turns out to be false. 

The non-invariance of the distorted Hamiltonian is especially surprising if we only consider the subset of canonical transformations that are induced by simple coordinate transformations of the generalised coordinates (i.e., point transformations). Often, the choice of generalised coordinates is treated only as a matter of convenience in formulating the equations of a model, and the effects of this choice on the accuracy of the solution are not considered. 

However, two symplectic numerical simulations of the same system simply written in different generalised coordinates will have different distorted Hamiltonians, and will thus yield different numerical results. This is the case that we will present in more detail in this section.


%

Switching to indicial notation, we use indices $i,k,l,m$ regarding the original coordinates and $\al,\be$ regarding the transformed coordinates. (Both sets of indices take values between $1, \ldots, d$, for a system with $d$ degrees of freedom.) A transformation of the generalised coordinates%
\begin{gather}
    \qb^{\al} = Q^{\al}(\qv);\quad \textnormal{with inverse } q^i = \left(\Qinv\right)^{i}(\qvb),
\end{gather}
induces a corresponding canonical transformation on the phase space \cite{goldstein2001classical} that transforms the generalised momenta as
\begin{gather}
    \pb_{\al} = p_{i} \left( \pdv{\left(\Qinv\right)^{i}}{\qb^{\al}} \right) \circ \Qv(\qv);\quad \textnormal{with inverse } p_{i} = \pb_{\al} \left( \pdv{Q^{\al}}{q^{i}} \right) \circ \Qvinv(\qvb).
\end{gather}

It is known that the distorted Hamiltonian $\Htilde$ corresponding to a symplectic numerical method differs from the original Hamiltonian by an $h$-dependent linear combination of elementary Hamiltonians \cite{sanzserna1991order}. For example, the symplectic Euler method\footnote{There exists an adjoint of the symplectic Euler method \eqref{eq:symE1}--\eqref{eq:symE2} presented here, which, somewhat confusingly, is often also called the symplectic Euler method \cite{leimkuhler2005simulating,hairer2006geometric}. The arguments presented in this article for the variant \eqref{eq:symE1}--\eqref{eq:symE2} transfer naturally to the adjoint variant.}
\begin{align}
    \pvjp &= \pvj - h H_q\leftf(\qvj,\pvjp\rightf)\label{eq:symE1}, \\
    \qvjp &= \qvj + h H_p\leftf(\qvj,\pvjp\rightf),
    \label{eq:symE2}
\end{align}
has the distorted Hamiltonian
\begin{gather}
    \Htilde = H - \frac{h}{2} H_p H_q + \Ordo{h^2}
    ,\label{eq:symEHtilde}
\end{gather}
where the inner product $H_p H_q$ is the only first-order elementary Hamiltonian.
Eq.\ \eqref{eq:symEHtilde} also shows that the symplectic Euler method is consistent and first-order accurate.

Naturally, for the transformed system, the distorted Hamiltonian corresponding to the method \eqref{eq:symE1}--\eqref{eq:symE2} becomes
\begin{gather}
    \Hbartilde = \Hbar - \frac{h}{2} \Hbar_{\pb} \Hbar_{\qb} + \Ordo{h^2}
    . \label{eq:symEHtildetrf}
\end{gather}

Now, we will show that
the first-order elementary Hamiltonian $H_p H_q$ is not invariant to a canonical transformation induced by an arbitrary point transformation. Then, we show how the distorted, transformed Hamiltonian differs from the distorted Hamiltonian in the case of the symplectic Euler method \eqref{eq:symE1}--\eqref{eq:symE2}.

For the transformed system, the two derivatives are
\begin{align}
    \left( \Hbar_{\pb} \right)^{\alpha}
    \equiv \pdv{\Hbar}{\pb_{\al}} &= \pdv{H}{p_i} \pdv{p_i}{\pb_{\al}} + \pdv{H}{q^k} \pdv{q^k}{\pb_{\al}}
    \nonumber \\ 
    &= \pdv{H}{p_i} \pdv{}{\pb_{\al}} \left( \pb_{\be} \pdv{Q^{\be}}{q^i} 
\right) 
    \nonumber \\
    &= \pdv{H}{p_i} \pdv{Q^{\al}}{q^i}
    ,\label{eq:Hbardp}
\end{align}
and
\begin{align}
    \left( \Hbar_{\qb} \right)_{\alpha}
    \equiv \pdv{\Hbar}{\qb^{\al}} &= \pdv{H}{p_i} \pdv{p_i}{\qb^{\al}} + \pdv{H}{q^k} \pdv{q^k}{\qb^{\al}}
    \nonumber \\ 
    &= \pdv{H}{p_i} \pdv{}{\qb^{\al}}
    \left(
    \pb_{\be} \pdv{Q^{\be}}{q^i} 
    \right) + \pdv{H}{q^k} \pdv{(\Qinv)^{k}}{\qb^{\al}}
    \nonumber \\
    &= \pdv{H}{p_i} \pb_{\be} \pdv[2]{Q^{\be}}{q^i}{q^l}
    \pdv{(\Qinv)^{l}}{\qb^{\al}}
    + \pdv{H}{q^k} \pdv{(\Qinv)^{k}}{\qb^{\al}}
    ,\label{eq:Hbardq}
\end{align}
from which
\begin{align}
    \Hbar_{\pb} \Hbar_{\qb}
    & \equiv
    \pdv{\Hbar}{\pb_{\al}} \pdv{\Hbar}{\qb^{\al}} 
    \nonumber \\
    &= \pdv{H}{p_i} \pdv{Q^{\al}}{q^i}
    \left(
    \pdv{H}{p_k} \pb_{\be} \pdv[2]{Q^{\be}}{q^k}{q^l}
    \pdv{(\Qinv)^{l}}{\qb^{\al}}
    +
    \pdv{H}{q^m} \pdv{(\Qinv)^{m}}{\qb^{\al}}
    \right)
    \nonumber \\
    &= \pdv{H}{p_i}\pdv{H}{p_k}\pb_{\be}\pdv[2]{Q^{\be}}{q^k}{q^i}
    + \pdv{H}{p_i}\pdv{H}{q^i},
\end{align}
or, in the original coordinates,
\begin{align}
    \Hbar_{\pb} \Hbar_{\qb}
    & \equiv
    \pdv{\Hbar}{\pb_{\al}} \pdv{\Hbar}{\qb^{\al}} 
    \nonumber \\
    &= 
    \underbrace{
        \pdv{H}{p_i}\pdv{H}{p_k} p_{l} \pdv{\left(\Qinv\right)^{l}}{\qb^{\be}}
        \pdv[2]{Q^{\be}}{q^k}{q^i}
    }_{\mathrel{=:}\!\; \deltaHpHq}
    + \pdv{H}{p_i}\pdv{H}{q^i}.\label{eq:Hbdpdq2}
\end{align}

Clearly, this shows that the elementary Hamiltonian $H_p H_q$ written in the original coordinates is equivalent to its counterpart $\Hbar_{\pb} \Hbar_{\qb}$ written in the transformed coordinate system, if the coordinate transformation $\Qv$ is affine. However, for a more general change of coordinates, the two expressions are not necessarily equivalent. (Note that the application of a non-affine $\Qv$ might still, in special cases%
\footnote{As an example, for the (non-affine) coordinate transformation
from polar to Cartesian coordinates, $\displaystyle \pdv{\left(\Qinv\right)^{l}}{\qb^{\be}} \pdv[2]{Q^{\be}}{q^k}{q^i} \not\equiv 0$. Yet, for the artificially constructed Hamiltonian $\displaystyle H(r,\pr,\pth) = \frac{1}{2} \left( p_r^2 + 2 \frac{\pth^2}{r^2} \right)$,
the invariance property
$\Hbar_{\pb} \Hbar_{\qb} = H_p H_q$ turns out to hold%
.},
result in 
$\deltaHpHq=0$%
.%
)

As the distorted Hamiltonian \eqref{eq:symEHtilde} of the symplectic Euler method contains this elementary Hamiltonian, we have shown that it is not invariant in general to changes of the coordinate system, i.e.\ $\Htilde \neq \Hbartilde$. Fig.~\ref{fig:commutative} also illustrates the non-invariance of the distorted Hamiltonian.

\begin{figure}[htbp]
    \centering
    \begin{tikzcd}[row sep=8ex, column sep=12ex]
        H
        \arrow[r, dashrightarrow, "\textnormal{discretization}"]
        \arrow[d, leftharpoonup, "\Qvinv",  shift left =0.20ex]
        \arrow[d, rightharpoondown, "\Qv"', shift right=0.20ex] 
        & \Htilde 
        \arrow[d, dash, dashed]
        \arrow[d, phantom, "\cross" marking]
        \\
        \Hbar \arrow[r, dashrightarrow, "\textnormal{discretization}"]
        & \Hbartilde
    \end{tikzcd}
    \caption{Relationship between the original Hamiltonian $H$, its transformed counterpart $\Hbar$, as well as their respective distorted versions $\Htilde$ and $\Hbartilde$ resulting from the application of a symplectic numerical method.}
    \label{fig:commutative}
\end{figure}

Actually, all elementary Hamiltonians are constructed from derivatives of $H$ of a certain order with respect to $\qv$ or $\pv$ (see \cite{sanzserna1991order,hairer1994backward,hairer2006geometric} for details, and esp.\ Section 5 of \cite{hairer1994backward} for the first few elementary Hamiltonians). This means that the above reasoning extends naturally to higher-order elementary Hamiltonians, as \eqref{eq:Hbardp}--\eqref{eq:Hbardq} show that while $H_p$ is invariant to the point transformation-induced canonical transformation, $H_q$ is not, and consequently, neither will its higher-order derivatives. As a result, the non-invariance discussed should be characteristic of other symplectic methods as well. However, a deeper investigation of higher-order methods is beyond the scope of the present article, as a rich structure emerges even for for first-order methods, which we will explore in more detail.

It is worth emphasizing here that since we have proven the general non-invariance of the distorted Hamiltonian for a subset of all possible canonical transformations, it follows immediately that a distorted Hamiltonian in general is not invariant with respect to all canonical transformations.

\subsection{Numerical demonstration of non-invariance}
To demonstrate the above results numerically, we return to the model of the planar elastic pendulum introduced in Section~\ref{sec:example_rowlands} (and detailed in \ref{app:elasticpendulum}). In the present example, we consider the numerical solution of the motion of the elastic pendulum formulated in a Cartesian and a polar coordinate system.

A natural question in light of the above would be to ask whether, for a given system and two coordinate systems, can one coordinate system be chosen over the other where the numerical results are ``better'' in a certain sense. We will see here and in Section~\ref{sec:firstintegrals} that the answer to this question depends fundamentally on the sense in which one method is considered ``better''. In this example, we consider the error in total energy as the measure of performance, as this is a natural benchmark of symplectic numerical methods.

We denote the mean $\ell_2$ error in total energy over time of the two numerical solutions in the two coordinate systems by $\epsxy$ and $\epsrphi$, respectively. As the total energy of the exact solution stays constant, this also measures the deviation of the corresponding distorted Hamiltonians from the original Hamiltonian. We have shown this deviation (i.e. $|H-\Htilde|$) should be expected to differ between coordinate systems if the transformation between them is not affine, which is indeed the case for the Cartesian and polar coordinate systems. One characteristic of this measure of performance is that it depends on the motion of the simulated systems, which in turn will depend on the initial condition as well.

\begin{figure}[htbp]
    \centering
    \includegraphics[trim=0 2mm 0 0 ,clip,width=0.65\textwidth]{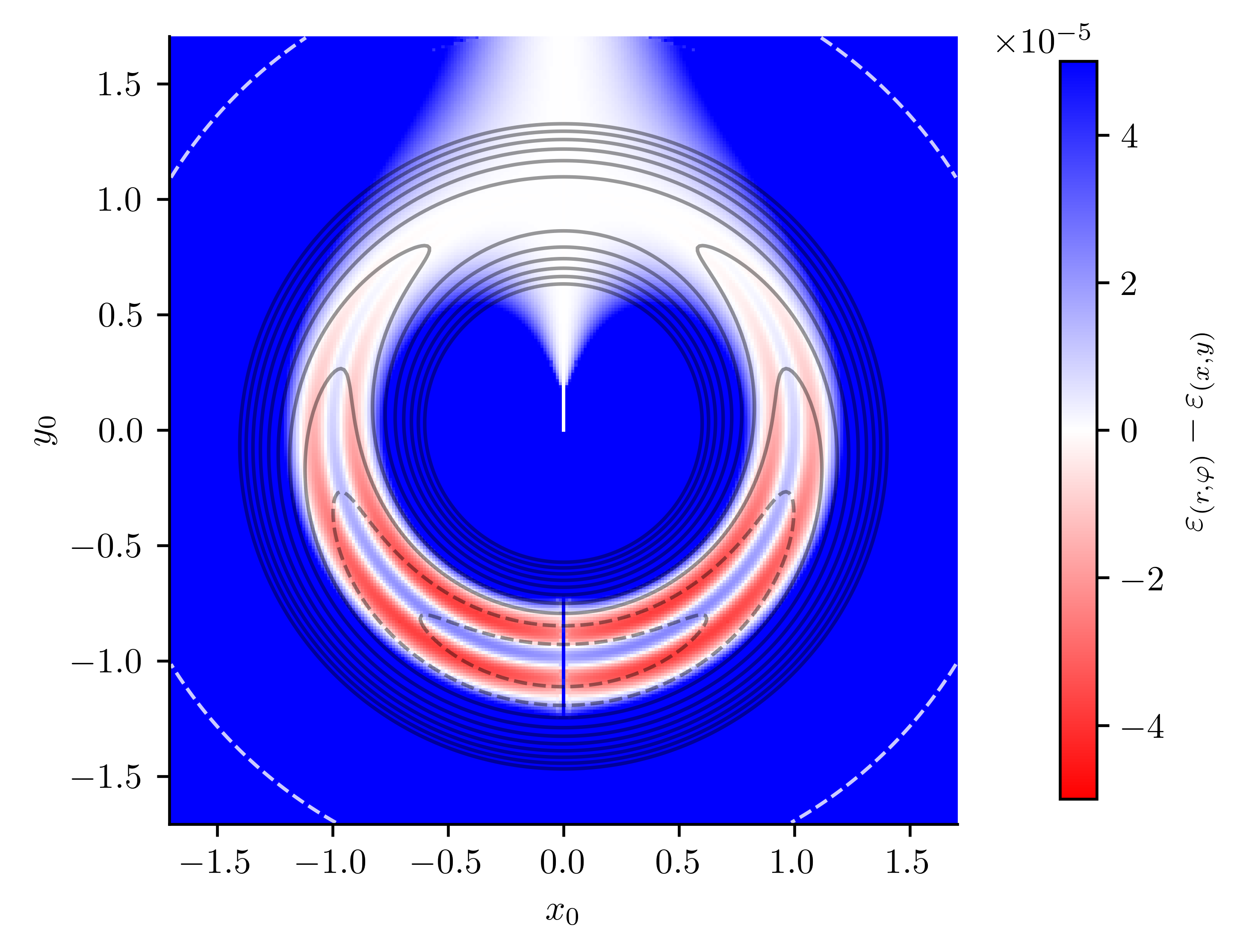}
        \caption{Comparison of error in total energy between simulations of the elastic pendulum performed in polar and Cartesian coordinates using the symplectic Euler method, with $H$-isolines (gray) and the boundary of stability for the polar version overlaid. (Parameters: $l{=}1,\, m{=}1,\, k{=}1,\, g{=}0.02,\, t_\mathrm{max}{=}50000,\, \Dt {=} 0.2.$)}
    \label{fig:error_coords_smooth_ic}
\end{figure}

For this reason, for all simulations, we use the initial conditions $p_x(0) = p_y(0) = 0$ for the generalized momenta, while the initial position is varied in the intervals $x(0) \in [-1.5,1.5]$, $y(0) \in [-1.5,1.5]$. Fig.~\ref{fig:error_coords_smooth_ic} shows the difference on the errors $\epsxy$ and $\epsrphi$ for long-term simulations as a function of the initial conditions. It is clear that though the energy-wise performance of the two methods differ, there is no clear winner: for some initial conditions (shown in blue), more accurate results are achieved in the Cartesian coordinates, and for others (shown in red) the polar coordinates are superior. This is true even if we only consider initial conditions corresponding the same total energy ($H$-isolines plotted as gray curves). We must conclude that from a total energy perspective, the best coordinate system, in general, can only be chosen when the initial condition is taken into account.

Additionally, it must also be noted that in this model, the simulations performed in a polar coordinate system diverge for initial conditions $p_x(0) = p_y(0) = 0$, $r(0) \geq 2 l + \frac{2 m g}{k} \cos\left( \varphi(0) \right)$, as the trajectory will pass through the singularity of the coordinate system at $r=0$. (The boundary of the region of unstable initial conditions is marked by a white dashed line in Fig.~\ref{fig:error_coords_smooth_ic}.) This highlights the fact that not only can the choice of the coordinate system have an effect on the accuracy of the numerical solution, but it can also influence its stability. However, a deeper analysis of the relationship between stability and the choice of coordinate system is beyond the scope of this article and is left for future work.

\subsection{Order compensation using coordinate transformation}\label{sec:ordercompensation}
Having de\-mon\-stra\-ted the non-invariance of the distorted Hamiltonian, we briefly investigate potential ways of exploiting it with the aim of increasing the order of accuracy of the numerical solution. 

First, let us return to \eqref{eq:symEHtildetrf}. If we know that $\Hbartilde$ can be different from $\Htilde$, we can see that one choice for a coordinate system that improves the convergence order of the method \eqref{eq:symE1}--\eqref{eq:symE2} and other symplectic methods would be a transformation that guarantees $\Hbar_{\qb} \equiv \Nv$. Component-wise, this is exactly the definition of cyclic coordinates \cite{goldstein2001classical} in Hamiltonian mechanics: thus, using the symplectic Euler method (or its adjoint variant) for a Hamiltonian problem rewritten in fully cyclic coordinates (action-angle variables) increases the order of accuracy. In fact, a slightly closer look reveals that using cyclic coordinates, the symplectic Euler method or any other numerical method is trivially an exact integrator of a system. Note that this line of thought has close parallels to that given in Section~\ref{sec:vtrf}.

However, it is well-known that most systems can not be transformed into a completely integrable form. Nevertheless, this is not the only option. For the symplectic Euler method, \eqref{eq:Hbdpdq2} gives an opportunity to determine a coordinate transformation $\Qv$ that ensures that $\Hbar_{\pb} \Hbar_{\qb} \equiv 0$ along the motion. Our experience is that the existence of such a coordinate transformation is not guaranteed, but the validity of the approach can be demonstrated through the example given in Appendix~\ref{app:harmonicoscillator}.

\section{Preservation of first integrals and change of generalised coordinates}\label{sec:firstintegrals}
So far, we have considered the effect of coordinate transformations on the distorted Hamiltonian and on the overall order of the method. However, other aspects of symplectic methods should also be considered that are connected to characteristic properties of Hamiltonian systems. One central aspect of numerical models used for simulating dynamical systems is the preservation of first integrals of the system that would be preserved along the motion. For a closed, autonomous system, the Hamiltonian itself is a first integral, but other symmetries of a system can give rise to additional first integrals (at most $n-1$ of them), e.g. the rotational symmetry of a system ensures the preservation of total angular momentum in a closed, autonomous system. The preservation of such a first integral in the numerical solution is not guaranteed: this depends on the numerical method used, as well as the coordinate system the solution is calculated in.

It is a well-known result in numerical analysis that among Runge--Kutta methods, symplectic ones preserve quadratic first integrals of the form
\begin{gather}
    \yv^{\Trp} \Cm \yv,
\end{gather}
(where $\Cm$ is a symmetric matrix), and vice versa: all Runge--Kutta methods that preserve such quadratic first integrals are symplectic \cite{hairer2008solving1,sanzserna1988runge}. An extension of this theorem exists for the more general class of B-series methods \cite{calvo1994canonical,chartier2006algebraic}, as well as regarding P-series methods, which describe partitioned numerical methods \cite{hairer2006geometric}. The symplectic Euler method used as an example throughout this article is a P-series method, thus it preserves quadratic first-integrals.

The above traditional formulation, however, sidesteps the important fact that no first integral is inherently quadratic: this depends on the variables (and the coordinate system) the model is formulated in, as the same first integral can be quadratic in one coordinate system and non-quadratic in others. Furthermore, a Hamiltonian system written using a cyclic coordinate will preserve the corresponding first integral trivially when solved using any numerical method, as it will be the generalized momentum corresponding to the aforementioned cyclic coordinate. At the same time, expressing the same system in a form without the cyclic coordinate results in equations where the first integral in question may not be preserved during numerical simulation.

Clearly, somewhere during the transformation of a Hamiltonian system from one set of variables to another, a numerical method can lose its ability to preserve one or more first integrals. Our aim in this section to investigate this with respect to canonical transformations induced by coordinate transformations for the symplectic Euler method.

\subsection{Preservation of first integrals in the symplectic Euler method}
Consider an autonomous system with a Hamiltonian $H$, expressed using generalised coordinates $q^i$ and generalised momenta $p_i$. Assume that the generalized momentum $p_1$ is a first integral that is preserved along the motion: this can be expressed using Poisson brackets as
\begin{gather}
    \{ H, p_1 \} = 0, \label{eq:p1preserved}
\end{gather}
which, when computed explicitly, gives
\begin{gather}
    \pdv{H}{q^1} = 0,
\end{gather}
showing that the corresponding generalized coordinate $q^1$ is indeed cyclic.

Now, the question we ask is the following: under which coordinate transformations (of the form $Q^{\al}(\qv)$) is $p_1$ preserved by a symplectic method at least up to order $r$? This condition can again be expressed using Poisson brackets, as
\begin{gather}
    \{ \Hbartilde, p_1 \} = 0 + \Ordo{h^r}.
\end{gather}

Here, we only treat the case of the symplectic Euler method, with $r=2$, as the method is originally first-order. Substituting the corresponding distorted transformed Hamiltonian \eqref{eq:symEHtildetrf} gives
\begin{gather}
    \left\{ \Hbar - \frac{h}{2} \Hbar_{\pb} \Hbar_{\qb} + \Ordo{h^2}, p_1 \right\} = 0 + \Ordo{h^2},
\end{gather}
which, after using the bilinearity of the Poisson brackets and \eqref{eq:p1preserved}, reduces to
\begin{gather}
    \left\{\Hbar_{\pb} \Hbar_{\qb}, p_1 \right\} = 0.
\end{gather}
Using the main result of the previous section, we can substitute the elementary Hamiltonian into the LHS as \eqref{eq:Hbdpdq2}, yielding
\begin{gather}
    \left\{ \pdv{H}{p_i}\pdv{H}{p_k} p_{l} \pdv{\left(\Qinv\right)^{l}}{\qb^{\be}} \pdv[2]{Q^{\be}}{q^k}{q^i}, p_1 \right\}
    + 
    \left\{ \pdv{H}{p_i}\pdv{H}{q^i}, p_1 
    \vphantom{\pdv{\left(\Qinv\right)^{l}}{\qb^{\be}}}
    \right\}
    = 0.
\end{gather}
Calculating the Poisson brackets explicitly gives
\begin{gather}
    \pdv{}{p_m}\left(\pdv{H}{p_i}\pdv{H}{p_k} p_{l} \pdv{\left(\Qinv\right)^{l}}{\qb^{\be}} \pdv[2]{Q^{\be}}{q^k}{q^i}\right) \pdv{p_1}{q^m}
    \mathrel{-} \nonumber \\
    \pdv{}{q^m}\left(\pdv{H}{p_i}\pdv{H}{p_k} p_{l} \pdv{\left(\Qinv\right)^{l}}{\qb^{\be}} \pdv[2]{Q^{\be}}{q^k}{q^i}\right) \pdv{p_1}{p_m}
    \mathrel{+} \nonumber \\
    \pdv{}{p_m}\left(\pdv{H}{p_i}\pdv{H}{q^i}\right) \pdv{p_1}{q^m}
    -
    \pdv{}{q^m}\left(\pdv{H}{p_i}\pdv{H}{q^i}\right) \pdv{p_1}{p_m} = 0,
\end{gather}
simplifying which yields
\begin{gather}
   \pdv{}{q^1}\left(\pdv{H}{p_i}\pdv{H}{p_k} p_{l} \pdv{\left(\Qinv\right)^{l}}{\qb^{\be}} \pdv[2]{Q^{\be}}{q^k}{q^i}
   \right)
    = 0.
\end{gather}
This, using the notation introduced in Section~\ref{sec:distortedhamiltonian}, can be also written as
\begin{gather}
    \pdv{\!\:\deltaHpHq}{q^1} = 0.\label{eq:firstintegral}
\end{gather}

Equation \eqref{eq:firstintegral} gives a necessary and sufficient condition for the preservation of the first integral $p_1$ in the transformed system up to second order. This means that \eqref{eq:firstintegral} is a necessary condition for the exact preservation of $p_1$: equivalently, if $\pdv{\!\:\deltaHpHq}{q^1} \neq 0$, then the first integral in question will not be preserved in the transformed system.

To demonstrate the application of this condition, we provide the following example.

\subsection{Example: free point mass in Cartesian and polar coordinates}
We consider the Hamiltonian system describing the motion of a single point mass in a plane, without the influence of any forces. Even in this elementary system, the transformation between the Cartesian and polar coordinates have a surprising effect on the preservation of first integrals when using the symplectic Euler method.

In Cartesian coordinates $\qv = (x,y)$, with corresponding generalized momenta $\pv = (p_x, p_y)$, the Hamiltonian can be written as
\begin{gather}
    H(\qv, \pv) = \frac{1}{2m} \left( p_x^2 + p_y^2 \right),
\end{gather}
where $m$ is the mass of the particle. Clearly, both $x$ and $y$ are cyclic coordinates here: $\partial H/\partial x = \partial H/\partial y \equiv 0$, meaning that both $p_x$ and $p_y$ are preserved. As these are in fact the components of the momentum vector, this corresponds to the preservation of linear momentum.

In polar coordinates $\qb = (r,\theta)$, however, the Hamiltonian of the same system can be shown to take the form of
\begin{gather}
    \Hbar(\qvb, \pvb) = \frac{1}{2m} \left( p_r^2 + \frac{p_{\theta}^2}{r^2} \right),
\end{gather}
with $\pv = (p_r, p_{\theta})$. Here, only $\theta$ is cyclic; and the preservation of $p_{\theta}$ is the preservation of the angular momentum.

Now, let us consider the coordinate transformations between these two coordinate systems written as
\begin{gather}
    \begin{pmatrix}
        r \\
        \theta
    \end{pmatrix}
    =
    \qvb = \Qv(\qv) = 
    \begin{pmatrix}
        \sqrt{x^2 + y^2} \\
        \atantwo(x,y)
    \end{pmatrix}
    \,
    \Leftrightarrow
    \,
    \begin{pmatrix}
        x \\
        y
    \end{pmatrix}
    =
    \qv = \Qvinv(\qvb) = 
    \begin{pmatrix}
        r \cos \theta \\
        r \sin \theta
    \end{pmatrix},\hspace*{-2em}
\end{gather}
where $\atantwo$ denotes the four-quadrant inverse tangent function.

\subsubsection{Preservation of linear momentum for a simulation in polar coordinates}
To get information about the preservation of the linear momentum components in a polar simulation, \eqref{eq:firstintegral} can be computed directly, with $x$ or $y$ as $q_1$. For $x$, the calculation yields
\begin{align}
    \pdv{\!\:\deltaHpHq}{x} \mathrel{=} & \frac{(\py x-\px y) \px^2 \left(y^3-3 x^2 y\right)}{\left(x^2+y^2\right)^3}
    \mathrel{+} \nonumber\\
    & \frac{(\py x-\px y) \bigleft(\px \py x \left(x^2-7 y^2\right)+2 \py^2 y (x-y) (x+y)\bigright)}{\left(x^2+y^2\right)^3}
\end{align}
and, for $y$,
\begin{align}
    \pdv{\!\:\deltaHpHq}{y} \mathrel{=} & -\frac{(\py x-\px y) \left(2 \px^2 x \left(y^2-x^2\right)\right)}{\left(x^2+y^2\right)^3},
    \mathrel{-} \nonumber\\
                                        & \frac{(\py x-\px y) \left(\px \py y \left(y^2-7 x^2\right)+\py^2 x \left(x^2-3 y^2\right)\right)}{\left(x^2+y^2\right)^3},
\end{align}
which are clearly not identically zero; thus we have shown that the linear momentum vector is not preserved in this case. 

To verify this statement, we can express e.g. $p_x$ in polar phase space coordinates as
\begin{gather}
    p_x = - \frac{p_{\theta}}{r} \sin \theta + p_r \cos \theta,
\end{gather}
and it is straightforward to show that after a time step, $p_x^{j+1} = p_x^{j} + \Ordo{h^2}$. (For $p_y$, we get a similar result.) Thus, it is clear that the condition \eqref{eq:firstintegral} indeed predicted the non-preservation of the components of the linear momentum.

Even though the linear momentum vector is not preserved by the symplectic Euler method in polar coordinates, the angular momentum is preserved trivially. However, in the reverse direction, the situation is different.

\subsubsection{Preservation of angular momentum for a simulation in Cartesian coordinates}
If the same system is simulated using the symplectic Euler method in Cartesian coordinates, the linear momentum vector is trivially preserved, and the question is the preservation of the angular momentum. Once again, \eqref{eq:firstintegral} can be used (with the role of the original and transformed coordinates reversed, and taking $\theta$ as $q_1$), yielding
\begin{gather}
    \pdv{\!\:\deltabarHpHq}{\theta} =
    0,\label{eq:firstintegralpolar}
\end{gather}
indicating that the angular momentum is preserved numerically at least up to the second order. In fact, this result is almost completely independent of the Hamiltonian: it can be shown that
\begin{gather}
    \pdv{}{\theta}\left(
        \pdv{Q^{\gamma}}{q^{i}} \pdv[2]{\left(\Qinv\right)^{i}}{\qb^{\al}}{\qb^{\be}}
    \right) \equiv 0,\quad \al,\be,\gamma = 1,2
\end{gather}
meaning that \eqref{eq:firstintegralpolar} is fulfilled trivially for any Hamiltonian for which $\theta$ is a cyclic coordinate.

Moreover, the angular momentum in Cartesian phase space coordinates can be expressed as
\begin{gather}
    p_{\theta} = x p_y- y p_x,
\end{gather}
which is a quadratic function of the phase space coordinates: thus, it is exactly preserved by the symplectic Euler method due to the aforementioned theorem on symplectic P-series.

In conclusion, we have demonstrated that even in a basic case such as this example, the choice of the coordinate system can greatly influence the preservation of first integrals of the system -- and the newly derived condition \eqref{eq:firstintegral} gives a sufficient condition to prove the non-preservation or the at least $\Ordo{h^2}$ preservation of a first integral.

\section{Conclusions}
We have given a general overview of the effect of variable transformations on ODEs and numerical methods used for solving them in general. We have shown through examples that exploiting these effects is possible, even though -- in general -- finding an optimal variable transformation can be of comparable difficulty to solving the original problem analytically. Then, we focused on symplectic methods used for simulating Hamiltonian systems, where an existing class of approaches, called processing methods, utilize canonical transformations close to the identity to achieve a higher rate of convergence.

After introducing the necessary framework, we extended the scope of the investigated transformations to include all canonical transformations that are induced by transformations defined for the generalised coordinates. For this class of coordinate transformations, we have shown that the distorted Hamiltonians corresponding to symplectic methods formulated in the original and the transformed variables are, in general, different -- even though the original Hamiltonian is invariant to canonical transformations by definition. As detailed in Section~\ref{sec:distortedhamiltonian}, we have shown this rigorously for the symplectic Euler method, while the arguments presented can be generalized to all symplectic methods, as the higher-order terms in their distorted Hamiltonians are functions of elementary Hamiltonians composed of further derivatives of the non-invariant terms discussed here.

This means that the accuracy of a symplectic method can depend on the choice of the coordinate system. We have demonstrated this result based on numerical simulations of the elastic pendulum, which also served as an example of stability being dependent on the coordinate system used. The effects of the choice of coordinates on stability also warrant further investigations in the future.

We have also investigated the possibility of exploiting the discrepancy in the distorted Hamiltonians to achieve a higher order of accuracy with the symplectic Euler method by choosing an optimal coordinate system, in Section~\ref{sec:ordercompensation}. Besides the trivial approach of using action-angle variables, we have also identified another approach, where the first-order term of the distorted Hamiltonian is eliminated along the motion. We were able to demonstrate the viability of this approach in Appendix~\ref{app:harmonicoscillator}. In the future, we would like to extend this approach and develop a systematic method for finding optimal coordinates.

We would like to take the opportunity to reflect briefly here on the choice of the total energy of the system along the motion over time as the measure of performance in Section~\ref{sec:distortedhamiltonian}. Many other aspects could be chosen to compare the behaviour of a symplectic numerical method in two different coordinate systems: e.g.\ we think that comparing the distorted vector fields directly would be a fruitful approach. However, such a comparison would need an appropriate metric defined on vector spaces of symplectic manifolds, which is still an open problem in symplectic geometry to the best of our knowledge. While there exits a relatively developed choice, Hofer's metric \cite{mcduff2017introduction}, it unfortunately can only be applied to compactly supported Hamiltonian maps, a condition that is not fulfilled by most physical systems -- with the exception of symplectic one-turn maps arising in the aperture calculation of particle accelerators \cite{erdelyi2001optimal}.

Besides the discrepancies in the distorted Hamiltonian itself, we have also investigated the preservation of first integrals in the symplectic Euler method in different coordinate systems, in Section~\ref{sec:firstintegrals}. We have derived a necessary condition \eqref{eq:firstintegral} for the exact preservation of a first integral corresponding to a cyclic coordinate, which depends both on the coordinate transformation considered and the Hamiltonian of the system involved. Naturally, this gives a sufficient condition for the non-preservation of said first integral. We have demonstrated said non-preservation condition through the example of the motion of a free point mass simulated in polar and Cartesian coordinates. This also served as an example of a case where the 
fulfilment of the aforementioned condition
is determined solely by the coordinate transformation, independently of the Hamiltonian.

As first integrals usually correspond to physical quantities that are important with respect to the practical goals of the simulation itself, we think it would be worthwhile to derive similar conditions for other symplectic numerical methods. Additionally, a sufficient and necessary condition for the exact preservation of a first integral corresponding to a cyclic coordinate may also be derived.

In this article, we focused on canonical transformations that are defined by a transformation of the generalised coordinates, as their choice can often be arbitrary or ad hoc, and a common attitude is to choose them based on the convenience of formulating the model, since they are often not presumed to affect the accuracy of the resulting simulation. As we have shown, however, the choice of coordinates can influence significantly both the overall accuracy of the results and the preservation of first integrals.

Additionally, as an extension of our present work, the effect of general canonical transformations on the accuracy of symplectic methods could also be studied, as they allow for even richer possibilities for finding an optimal coordinate system.

\section*{Acknowledgement}
This work has been supported by the EK\"{O}P-24-3-BME-247 grant of the National Research, Development and Innovation Office under the University Researcher Scholarship Program (EK\"{O}P)

We also thank Gabriella Svantnern\'{e} Sebesty\'{e}n and Istv\'{a}n Farag\'{o} for their comments.

\section*{Author contributions (CRediT)}
D. M. Tak\'{a}cs: Conceptualization, Methodology, Software, Formal Analysis, Writing -- Original Draft, Visualization. T. F\"{u}l\"{o}p: Formal analysis, Writing -- Review \& Editing, Supervision.


\cleardoublepage

\appendix

\section{The planar elastic pendulum}\label{app:elasticpendulum}

Here, for reference, we give the expressions for the mechanical model of the planar elastic pendulum.

\begin{figure}[htbp]
    \centering
    \includegraphics[width=0.35\textwidth]{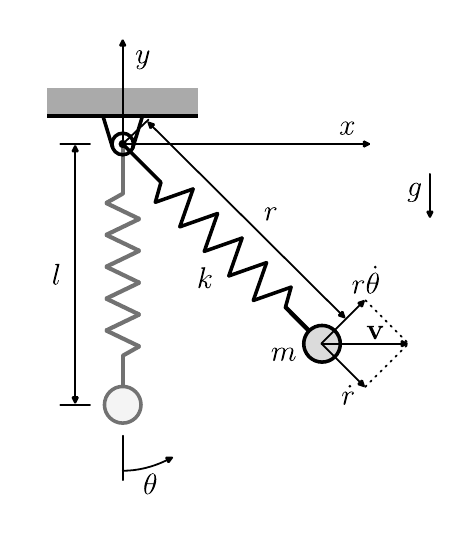}
    \caption{Illustration of the mechanical model of the planar elastic pendulum.}
    \label{fig:inga}
\end{figure}

 Fig.~\ref{fig:inga} illustrates the mechanical model of the planar elastic pendulum with point mass $m$ stiffness $k$ and relaxed length $l$, subject to gravitational acceleration $g$. Two coordinate systems are chosen for describing the motion, a Cartesian coordinate system with coordinates $(x, y)$, and a polar coordinate system with coordinates $(r, \theta)$.

In polar coordinates, the Lagrangian $\Lcal$ of the planar elastic pendulum can be expressed as
\begin{align}
    \Lcal(r, \theta, \dot{r}, \dot{\theta}) &= \Tcal(r, \dot{r}, \dot{\theta}) - \Ucal(r, \theta) \nonumber\\
                                            &= \frac{1}{2} m \left( \dot{r} + r^2 \dot{\varphi}^2 \right) - \frac{1}{2} k \left( r-l \right)^2 + m g r \cos(\theta),
\end{align}
which allows for the expression of the generalized momenta as
\begin{gather}
    \pr = \frac{\partial \Lcal}{\partial \dot{r}} = m \dot{r},
    \qquad
    \pth = \frac{\partial \Lcal}{\partial \dot{\theta}} = m r^2 \dot{\theta}.
\end{gather}
Consequently, the Hamiltonian in polar coordinates can be expressed as
\begin{align}
    \Hcal(r, \theta, \pr, \pth) &= \Tcal(r, \dot{r}, \dot{\theta}) + \Ucal(r, \theta) \nonumber \\
                                &= \frac{1}{2m}\left( \pr^2 + \frac{\pth^2}{r^2} \right) + \frac{1}{2} k \left( r-l \right)^2 - m g r \cos(\theta).\label{eq:peppol}
\end{align}
Following similar steps, the canonical momenta in a Cartesian coordinate system can be expressed as $\px=m\dot{x}$ and $\py=m\dot{y}$, with Hamiltonian
\begin{gather}
    \Hcal(x, y, \px, \py) = \frac{1}{2m}\left( \px^2  + \py^2 \right) + \frac{1}{2} k \left( \sqrt{x^2+y^2} - l \right)^2 - m g y\label{eq:pepdesc}.
\end{gather}

Note that the Rowlands method can only be applied for the equations \eqref{eq:pepdesc} formulated in the Cartesian coordinate system, due to the condition \eqref{eq:rowlandscondition}.

\section{Order compensation of the symplectic Euler method for the harmonic oscillator using an optimal coordinate transformation}\label{app:harmonicoscillator}

Here, we give an example of applying the basic idea of order compensation from Section \ref{sec:ordercompensation} on the harmonic oscillator, simulated using the symplectic Euler method.

In a dimensionless form, the Hamiltonian of the harmonic oscillator can be expressed naturally as
\begin{gather}
    H(q,p) = \frac{1}{2} p^2 + \frac{1}{2} q^2,\label{eq:harmoscH}
\end{gather}
yielding the canonical equations of motion
\begin{align}
    \dot{p} &= -q, \\
    \dot{q} &= p,
\end{align}
Without loss of generality, we can choose the initial condition $q(0) = 1,\,p(0) =0$, as we know the solution is periodic, and the amplitude has been used for deriving the dimensionless form \eqref{eq:harmoscH}.

We are looking to eliminate first-order errors along the motion, using an appropriately chosen coordinate transform. The condition for this is
\begin{gather}
    \Hbar - \Hbartilde \equiv \Ordo{h^2}
\end{gather}
which, based on \eqref{eq:symEHtildetrf}, can be simplified to
\begin{gather}
    \Hbar_{\pb} \Hbar_{\qb} \equiv 0.
\end{gather}
Using \eqref{eq:Hbdpdq2} and \eqref{eq:harmoscH}, in the original coordinate system this condition becomes
\begin{gather}
    p \left( p^2 \left(\pdv{Q}{q}\right)^{-1} \pdv[2]{Q}{q} + q \right) \equiv 0.
\end{gather}
A sufficient condition for this expression to hold is for the expression inside the parentheses to vanish. Using \eqref{eq:harmoscH} and the initial condition, the above condition along the solution is fulfilled by solutions of the differential equation
\begin{gather}
     \left( 1 - q^2 \right) \left(\pdv{Q}{q}\right)^{-1} \pdv[2]{Q}{q} + q = 0.
\end{gather}
For $-1 \leq q \leq +1$, with an appropriate selection of the the integration constants, this has a solution\footnote{It is worth noting that the RHS of \eqref{eq:Qexact} can be transformed into the radial Kepler equation for unit eccentricity \cite{walters2018astronautics} by the substitution $q = \sin(u)$.}
\begin{gather}
    Q(q) = \frac{2}{\pi} \left( q \sqrt{1-q^2} + \arcsin\leftf(q\rightf) \right)\label{eq:Qexact}
\end{gather}
as an optimal coordinate transform. However, the numerical solution will take values outside this interval, thus this coordinate transformation can not be used as is. Instead, we can take a modification that is $\Ordo{h^2}$-close, e.g.
\begin{gather}
    Q(q) = \frac{2}{\pi} \left( \hat{q} \sqrt{1-\hat{q}^2} + \arcsin\leftf(\hat{q}\rightf) \right),\,\textnormal{where }\hat{q} = \frac{q}{1+k h^2},\label{eq:Qhkozel}
\end{gather}
with $k>0$ being a parameter that should be chosen such that it guarantees $|\hat{q}| \leq 1$ for the numerical solution. It can be shown that the difference between \eqref{eq:Qexact} and \eqref{eq:Qhkozel} is indeed $\Ordo{h^2}$, and $k=2$ is an appropriate choice.

\begin{figure}[hbt]
    \centering
    \includegraphics[width=0.4\textwidth]{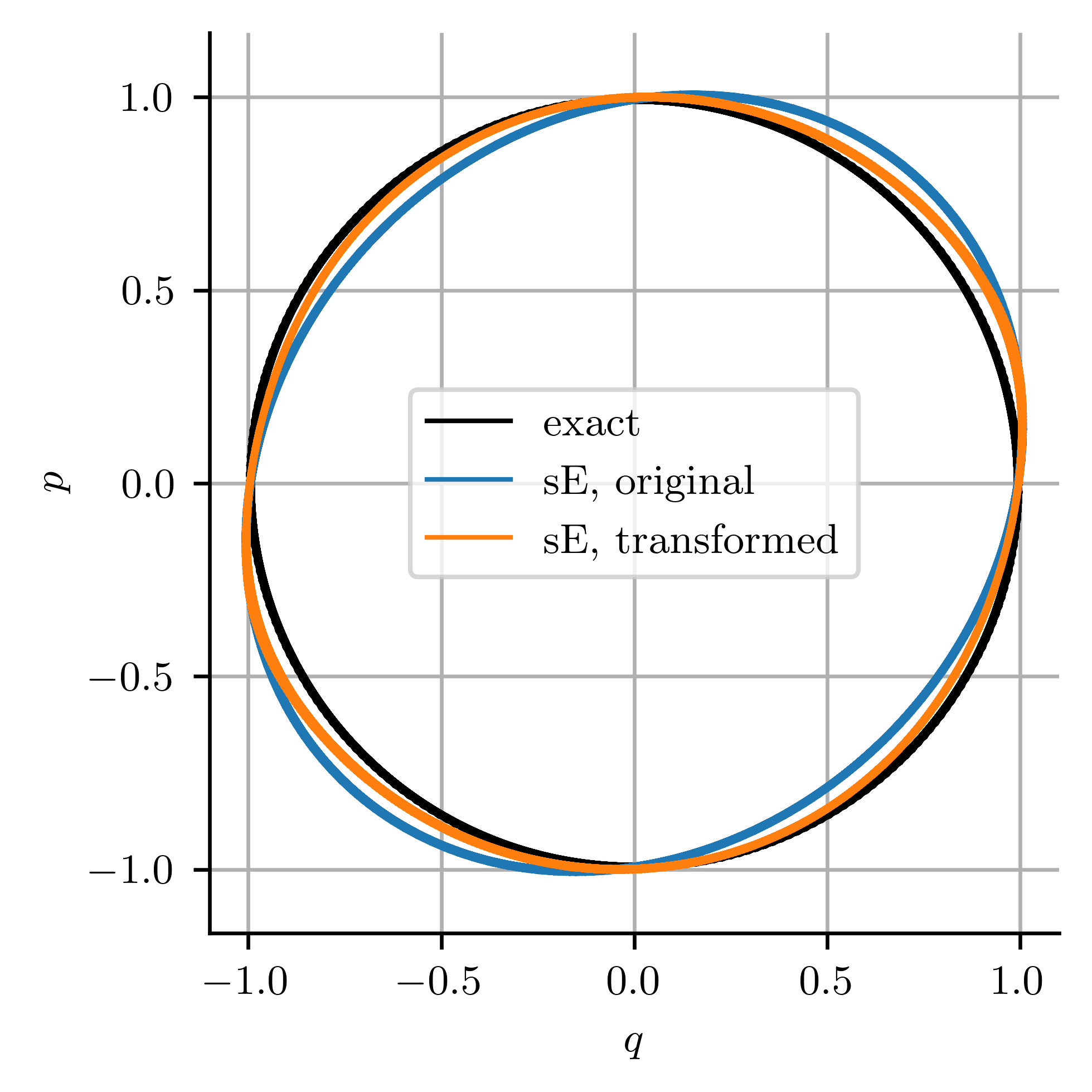}
    \caption{Trajectory of symplectic Euler simulation of the harmonic oscillator, with and without using the compensating coordinate transformation. (Parameters: $h=0.3$, $k=2$, number of steps: 1000.)}
    \label{fig:harmonic_oscillator_transformed_phasespace}
\end{figure}

\begin{figure}[hbt]
    \centering
    \includegraphics[width=0.4\textwidth]{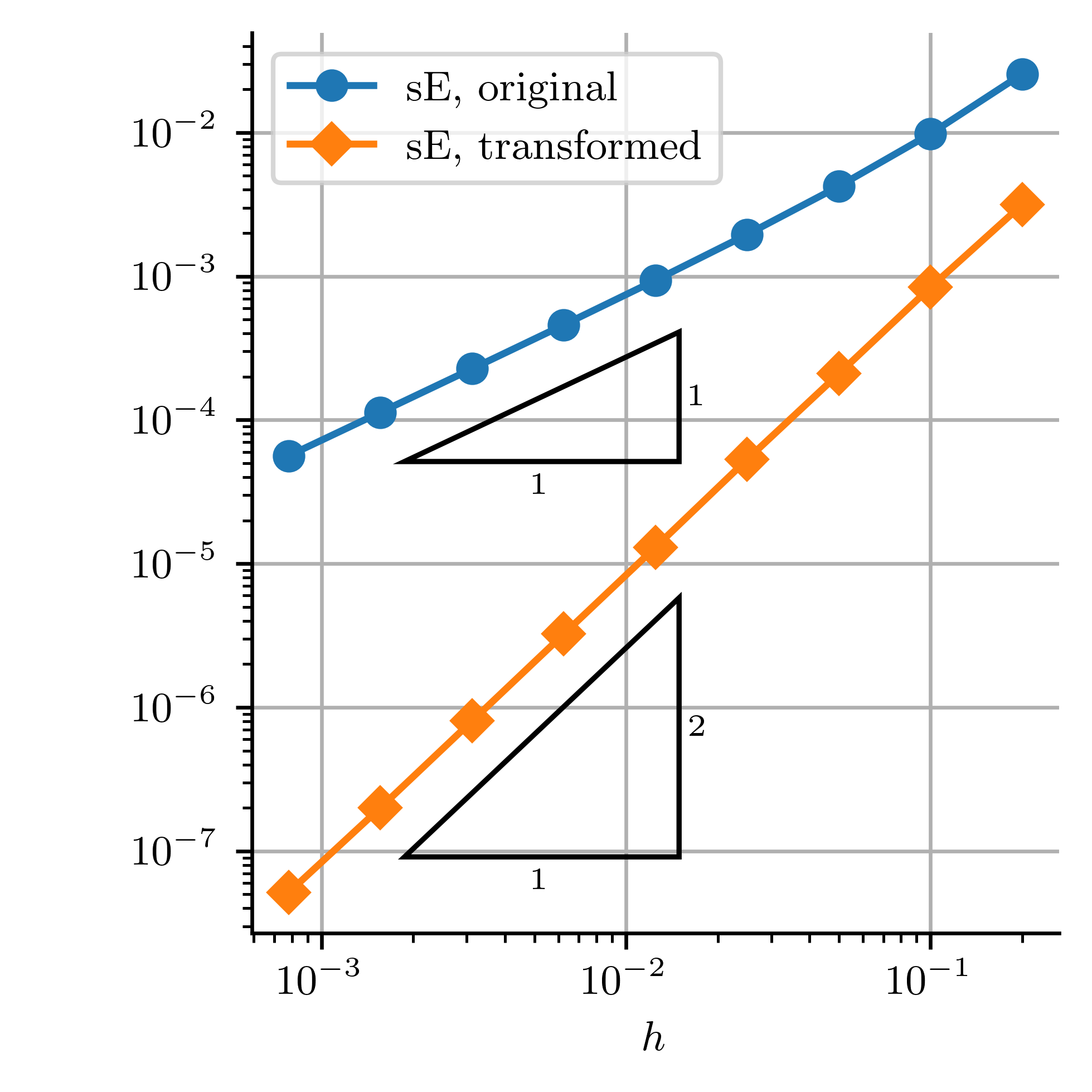}
    \caption{Convergence of the error in the Hamiltonian for symplectic Euler simulations of the harmonic oscillator with and without the compensating coordinate transformation.}
    \label{fig:harmonic_oscillator_transformed_convergence}
\end{figure}

Fig.~\ref{fig:harmonic_oscillator_transformed_phasespace} compares the performance of the symplectic Euler method with and without applying \eqref{eq:Qhkozel}, showing that the former is indeed closer to the exact solution. Fig.~\ref{fig:harmonic_oscillator_transformed_convergence} demonstrates that performing the simulation in the coordinate system defined by \eqref{eq:Qhkozel} indeed eliminates the first-order error from the total energy of the numerical solution, leading to a second-order convergence.

\printbibliography
\end{document}